\documentclass[a4paper,11pt,reqno]{amsart}
\usepackage{graphicx,latexsym,amssymb,xypic}

\newcommand{\sym}{\mathcal{S}} 
\newcommand{\N}{{\mathbb N}}

\newcommand{\ns}{\hspace{-0.2ex}} %Very thin negative space.

\theoremstyle{plain} 
\newtheorem{theorem}{Theorem}
\newtheorem*{theorem*}{Theorem}
\newtheorem{corollary}[theorem]{Corollary}
\newtheorem*{corollary*}{Corollary}
\newtheorem{porism}[theorem]{Porism} 
\newtheorem*{porism*}{Porism}
\newtheorem{lemma}[theorem]{Lemma} 
\newtheorem*{lemma*}{Lemma}
\newtheorem{proposition}[theorem]{Proposition}
\newtheorem*{proposition*}{Proposition}

\theoremstyle{definition} 
\newtheorem{definition}[theorem]{Definition}
\newtheorem*{definition*}{Definition}
\newtheorem{example}[theorem]{Example} 
\newtheorem*{example*}{Example}

\theoremstyle{remark} 

\newtheorem*{remark*}{Remark}

\newcommand{\wpi}{\widehat\pi} 
\newcommand{\emptyword}{{\epsilon}}

\DeclareMathOperator{\proj}{\mathrm{proj}}
\DeclareMathOperator{\des}{des}

\newcommand{\mono}{\mathcal{M}}

\def\dd{\hspace{1pt}\mbox{\textup{\text{-}\makebox[0pt]{\text{-}}}}\hspace{2pt}}

\def\ab{\ensuremath{ab}}

\def\axbc{\ensuremath{a{\dd}bc}}
\def\axcb{\ensuremath{a{\dd}cb}}
\def\bxac{\ensuremath{b{\dd}ac}}
\def\bxca{\ensuremath{b{\dd}ca}}
\def\cxab{\ensuremath{c{\dd}ab}}
\def\cxba{\ensuremath{c{\dd}ba}}

\def\abxc{\ensuremath{ab{\dd}c}}
\def\acxb{\ensuremath{ac{\dd}b}}
\def\baxc{\ensuremath{ba{\dd}c}}
\def\bcxa{\ensuremath{bc{\dd}a}}
\def\caxb{\ensuremath{ca{\dd}b}}
\def\cbxa{\ensuremath{cb{\dd}a}}

\def\axcxb{\ensuremath{a{\dd}c{\dd}b}}
\def\bxaxc{\ensuremath{b{\dd}a{\dd}c}}

\newcommand{\size}{1.1ex} 
\newcommand{\xxx}{\makebox[\size]{ }}

\newcommand{\CIR}{\makebox[\size]{\ensuremath{\circ}}}

\newcommand{\DAS}{\makebox[\size]{\ensuremath{-}}}

\newcommand{\ONE}{\makebox[\size]{1}}
\newcommand{\TWO}{\makebox[\size]{2}}
\newcommand{\THR}{\makebox[\size]{3}}
\newcommand{\FOU}{\makebox[\size]{4}}
\newcommand{\FIV}{\makebox[\size]{5}}
\newcommand{\SIX}{\makebox[\size]{6}}
\newcommand{\SEV}{\makebox[\size]{7}}
\newcommand{\EIG}{\makebox[\size]{8}}
\newcommand{\NIN}{\makebox[\size]{9}}
\newcommand{\Axx}{\makebox[\size]{10}}
\newcommand{\Bxx}{\makebox[\size]{11}}
\newcommand{\Cxx}{\makebox[\size]{12}}
\newcommand{\Dxx}{\makebox[\size]{13}}

\newcommand{\FIGmono}{
\begin{smallmatrix}
\xxx\xxx\xxx\xxx\xxx\xxx\xxx\xxx\xxx\xxx\xxx\xxx\xxx\xxx\xxx\xxx\xxx\xxx\xxx\xxx\xxx\xxx\xxx\xxx\CIR\xxx\xxx\CIR\DAS\DAS\CIR\DAS\DAS\DAS\DAS\DAS\CIR\\
\xxx\xxx\xxx\xxx\xxx\xxx\xxx\xxx\xxx\CIR\DAS\DAS\DAS\DAS\DAS\CIR\DAS\DAS\DAS\DAS\DAS\DAS\DAS\DAS\DAS\DAS\DAS\DAS\DAS\DAS\DAS\DAS\DAS\CIR\xxx\xxx\xxx\\
\xxx\xxx\xxx\xxx\xxx\xxx\CIR\DAS\DAS\DAS\DAS\DAS\DAS\DAS\DAS\DAS\DAS\DAS\DAS\DAS\DAS\CIR\xxx\xxx\xxx\xxx\xxx\xxx\xxx\xxx\xxx\xxx\xxx\xxx\xxx\xxx\xxx\\
\CIR\DAS\DAS\CIR\DAS\DAS\DAS\DAS\DAS\DAS\DAS\DAS\CIR\DAS\DAS\DAS\DAS\DAS\CIR\xxx\xxx\xxx\xxx\xxx\xxx\xxx\xxx\xxx\xxx\xxx\xxx\xxx\xxx\xxx\xxx\xxx\xxx\\
\ONE\xxx\xxx\TWO\xxx\xxx\THR\xxx\xxx\FOU\xxx\xxx\FIV\xxx\xxx\SIX\xxx\xxx\SEV\xxx\xxx\EIG\xxx\xxx\NIN\xxx\xxx\Axx\xxx\xxx\Bxx\xxx\xxx\Cxx\xxx\xxx\Dxx\\
\end{smallmatrix}
}

\newcommand{\FIGnop}{
\begin{smallmatrix}
\xxx\xxx\xxx\xxx\xxx\xxx\xxx\xxx\xxx\CIR\DAS\DAS\DAS\DAS\DAS\CIR\DAS\DAS\CIR\xxx\xxx\xxx\xxx\xxx\xxx\xxx\xxx\xxx\xxx\xxx\xxx\xxx\xxx\xxx\xxx\xxx\xxx\\
\xxx\xxx\xxx\xxx\xxx\xxx\CIR\DAS\DAS\DAS\DAS\DAS\DAS\DAS\DAS\DAS\DAS\DAS\DAS\DAS\DAS\CIR\xxx\xxx\CIR\xxx\xxx\CIR\DAS\DAS\CIR\DAS\DAS\CIR\xxx\xxx\xxx\\
\CIR\DAS\DAS\CIR\DAS\DAS\DAS\DAS\DAS\DAS\DAS\DAS\CIR\DAS\DAS\DAS\DAS\DAS\DAS\DAS\DAS\DAS\DAS\DAS\DAS\DAS\DAS\DAS\DAS\DAS\DAS\DAS\DAS\DAS\DAS\DAS\CIR\\
\ONE\xxx\xxx\TWO\xxx\xxx\THR\xxx\xxx\FOU\xxx\xxx\FIV\xxx\xxx\SIX\xxx\xxx\SEV\xxx\xxx\EIG\xxx\xxx\NIN\xxx\xxx\Axx\xxx\xxx\Bxx\xxx\xxx\Cxx\xxx\xxx\Dxx\\
\end{smallmatrix}
}

\begin{document}

\title{Generalised pattern avoidance}

\author{Anders Claesson}
\address{Matematik\\
  Chalmers tekniska h\"ogskola och G\"oteborgs universitet\\
  S-412 96 G\"oteborg, Sweden} \email{claesson@math.chalmers.se}
%\thanks{thanks} 
%\keywords{keywords}
%\subjclass{Primary: subject; Secondary: subject}
\date{\today}

\begin{abstract}
  Recently, Babson and Steingr\'{\i}msson have introduced generalised
  permutation patterns that allow the requirement that two adjacent
  letters in a pattern must be adjacent in the permutation. We
  consider pattern avoidance for such patterns, and give a complete
  solution for the number of permutations avoiding any single pattern
  of length three with exactly one adjacent pair of letters. We also
  give some results for the number of permutations avoiding two
  different patterns. Relations are exhibited to several well studied
  combinatorial structures, such as set partitions, Dyck paths,
  Motzkin paths, and involutions. Furthermore, a new class of set
  partitions, called monotone partitions, is defined and shown to be
  in one-to-one correspondence with non-overlapping partitions.
\end{abstract}

\maketitle\thispagestyle{empty}

\section{Introduction}

In the last decade a wealth of articles has been written on the
subject of pattern avoidance, also known as the study of ``restricted
permutations'' and ``permutations with forbidden subsequences''.
Classically, a pattern is a permutation $\sigma\in\sym_k$, and a
permutation $\pi\in\sym_n$ avoids $\sigma$ if there is no subsequence
in $\pi$ whose letters are in the same relative order as the letters
of $\sigma$. For example, $\pi\in\sym_n$ avoids $132$ if there is no
$1\leq i\leq j\leq k\leq n$ such that $\pi(i)\leq\pi(k)\leq\pi(j)$.
In \cite{Kn73v1} Knuth established that for all $\sigma\in\sym_3$,
the number of permutations in $\sym_n$ avoiding $\sigma$ equals the
$n$th Catalan number, $C_n = \frac{1}{1+n}\binom{2n}{n}$. One may also
consider permutations that are required to avoid several patterns. In
\cite{SiSc85} Simion and Schmidt gave a complete solution for
permutations avoiding any set of patterns of length three. Even
patterns of length greater than three have been considered. For
instance, West showed in \cite{We95} that permutations avoiding
both $3142$ and $2413$ are enumerated by the Shr\"oder numbers, $S_n =
\sum_{i=0}^n\binom {2n-i} i C_{n-i}$.

In \cite{BaSt00} Babson and Steingr\'{\i}msson introduced generalised
permutation patterns that allow the requirement that two adjacent
letters in a pattern must be adjacent in the permutation.  The
motivation for Babson and Steingr\'{\i}msson in introducing these
patterns was the study of Mahonian statistics, and they showed that
essentially all Mahonian permutation statistics in the literature can
be written as linear combinations of such patterns. An example of a
generalised pattern is $(\axcb)$. An $(\axcb)$-subword of a
permutation $\pi = a_1 a_2 \cdots a_n$ is a subword $a_i a_j a_{j+1}$,
$(i<j)$, such that $a_i<a_{j+1}<a_j$. More generally, a pattern $p$ is
a word over the alphabet $a<b<c<d\cdots$ where two adjacent letters
may or may not be separated by a dash. The absence of a dash between
two adjacent letters in a $p$ indicates that the corresponding letters
in a $p$-subword of a permutation must be adjacent. Also, the ordering
of the letters in the $p$-subword must match the ordering of the
letters in the pattern. This definition, as well as any other
definition in the introduction, will be stated rigorously in
Section~\ref{prel}. All classical patterns are generalised patterns
where each pair of adjacent letters is separated by a dash. For
example, the generalised pattern equivalent to $132$ is $(\axcxb)$.

We extend the notion of pattern avoidance by defining that a
permutation avoids a (generalised) pattern $p$ if it does not contain
any $p$-subwords. We show that this is a fruitful extension, by
establishing connections to other well known combinatorial structures,
not previously shown to be related to pattern avoidance.  The main
results are given below.
\vspace{.3ex}
$$
\vspace{.3ex}
\begin{array}{l|l|l}
  P & |\sym_n(P)| & Description\\ \hline
  \axbc                 &  B_n   & \text{ Partitions of } [n] \\
  \axcb                 &  B_n   & \text{ Partitions of } [n] \\
  \bxac                 &  C_n   & \text{ Dyck paths of length } 2n \\
  \axbc,\,\abxc         &  B^*_n & \text{ Non-overlapping partitions of } [n] \\
  \axbc,\,\axcb         &  I_n   & \text{ Involutions in } \sym_n  \\
  \axbc,\,\acxb         &  M_n   & \text{ Motzkin paths of length } n \\
\end{array}
$$
Here $\sym_n(P) = \{ \pi\in\sym_n : \pi \text{ avoids } p \text{ for
  all } p\in P \}$, and $[n]=\{1,2,\ldots,n\}$. When proving that
$|\sym_n(\axbc,\,\abxc)| = B^*_n$ (the $n$th Bessel number), we first
prove that there is a one-to-one correspondence between
$\{\axbc,\abxc\}$-avoiding permutations and \emph{monotone
  partitions}. A partition is monotone if its non-singleton blocks can
be written in increasing order of their least element and increasing
order of their greatest element, simultaneously. This new class of
partitions is then shown to be in one-to-one correspondence with
non-overlapping partitions.

\section{Preliminaries}\label{prel}

By an \emph{alphabet} $X$ we mean a non-empty set. An element of $X$
is called a \emph{letter}. A \emph{word} over $X$ is a finite sequence
of letters from $X$. We consider also the \emph{empty word}, that is,
the word with no letters; it is denoted by $\emptyword$. Let
$x=x_1x_2\cdots x_n$ be a word over $X$. We call $|x|:=n$ the
\emph{length} of $x$. A \emph{subword} of $x$ is a word
$v=x_{i_1}x_{i_2}\cdots x_{i_k}$, where $1 \leq i_1<i_2<\cdots<i_k
\leq n$.  A \emph{segment} of $x$ is a word $v=x_{i}x_{i+1}\cdots
x_{i+k}$. If $X$ and $Y$ are two linearly ordered alphabets, then two
words $x=x_1x_2 \cdots x_n$ and $y = y_1 y_2 \cdots y_n$ over $X$ and
$Y$, respectively, are said to be \emph{order equivalent} if $x_i<x_j$
precisely when $y_i < y_j$.

Let $X = A \cup \{\dd\}$ where $A$ is a linearly ordered alphabet. For
each word $x$ let $\bar x$ be the word obtained from $x$ by deleting
all dashes in $x$. A word $p$ over $X$ is called a \emph{pattern} if
it contains no two consecutive dashes and $\bar{p}$ has no repeated
letters. By slight abuse of terminology we refer to the \emph{length
  of a pattern} $p$ as the length of $\bar{p}$. Two patterns $p$ and
$q$ of equal length are said to be \emph{dash equivalent} if the $i$th
letter in $p$ is a dash precisely when the $i$th letter in $q$ is a
dash. If $p$ and $q$ are dash and order equivalent, then $p$ and
$q$ are \emph{equivalent}. In what follows a pattern will usually be
taken to be over the alphabet $\{a,b,c,d,\ldots\}\cup\{\dd\}$ where
$\{a,b,c,d,\ldots\}$ is ordered so that $a<b<c<d<\cdots$.

Let $[n]:=\{1,2,\ldots,n\}$ (so $[0]=\emptyset$). A \emph{permutation}
of $[n]$ is bijection from $[n]$ to $[n]$. Let $\sym_n$ be the set of
permutations of $[n]$.  We shall usually think of a permutation $\pi$
as the word $\pi(1)\pi(2)\cdots\pi(n)$ over the alphabet $[n]$.  In
particular, $\sym_0 = \{\emptyword\}$, since there is only one
bijection from $\emptyset$ to $\emptyset$, the empty map.  We say that
a subword $\sigma$ of $\pi$ is a $p$-\emph{subword} if by replacing
(possibly empty) segments of $\pi$ with dashes we can obtain a pattern
$q$ equivalent to $p$ such that $\bar q = \sigma$.  However, all
patterns that we will consider will have a dash at the beginning and
one at the end.  For convenience, we therefore leave them out.  For
example, $(\axbc)$ is a pattern, and the permutation $491273865$
contains three $(\axbc)$-subwords, namely $127$, $138$, and $238$. A
permutation is said to be \emph{$p$-avoiding} if it does not contain
any $p$-subwords. Define $\sym_n(p)$ to be the set of $p$-avoiding
permutations in $\sym_n$ and, more generally, $\sym_n(A) =
\bigcap_{p\in A} \sym_n(p)$.

We may think of a pattern $p$ as a permutation statistic, that is,
define $p\,\pi$ as the number of $p$-subwords in $\pi$, thus regarding
$p$ as a function from $\sym_n$ to $\N$. For example,
$(\axbc)\,491273865 = 3$. In particular, $\pi$ is $p$-avoiding if and
only if $p\,\pi = 0$. We say that two permutation statistics
$\mathrm{stat}$ and $\mathrm{stat}'$ are \emph{equidistributed} over
$A\subseteq\sym_n$, if
$$
\sum_{\pi\in A}x^{\mathrm{stat}\, \pi} = \sum_{\pi\in
  A}x^{\mathrm{stat}'\, \pi}.
$$
In particular, this definition applies to patterns. 

Let $\pi=a_1 a_2 \cdots a_n \in \sym_n$. An $i$ such that
$a_i>a_{i+1}$ is called a \emph{descent} in $\pi$. We denote by $\des
\pi$ the number of descents in $\pi$.  Observe that $\des$ can be
defined as the pattern $(ba)$, that is, $\des\pi = (ba)\pi$.  A
\emph{left-to-right minimum} of $\pi$ is an element $a_i$ such that
$a_i < a_j$ for every $j<i$. The number of left-to-right minima is a
permutation statistic.  Analogously we also define \emph{left-to-right
  maximum}, \emph{right-to-left minimum}, and \emph{right-to-left
  maximum}.

In this paper we will relate permutations avoiding a given set of
patterns to other better known combinatorial structures. Here follows
a brief description of these structures.

\subsection*{Set partitions}

A \emph{partition} of a set $S$ is a family, $\pi = \{A_1 , A_2,
\ldots, A_k \}$, of pairwise disjoint non-empty subsets of $S$ such
that $S=\cup_i A_i$. We call $A_i$ a \emph{block} of $\pi$.  The total
number of partitions of $[n]$ is called a \emph{Bell number} and is
denoted $B_n$. For reference, the first few Bell numbers are
$$
1,1,2,5,15,52,203,877,4140,21147,115975,678570,4213597. $$
Let
$S(n,k)$ be the number of partitions of $[n]$ into $k$ blocks; these
numbers are called the \emph{Stirling numbers of the second kind}.

\subsection*{Non-overlapping partitions}

Two blocks $A$ and $B$ of a partition $\pi$ \emph{overlap} if
$$\min A < \min B < \max A < \max B.$$
A partition is
\emph{non-overlapping} if no pairs of blocks overlap. 
Thus 
$$\pi = \{\{ 1,2,5,13\},\{3,8\},\{4,6,7\},\{9\},\{10,11,12\}\}$$
is non-overlapping. A pictorial representation of 
$\pi$ is
$$\pi \,=\, \FIGnop .$$
Let $B^*_n$ be
the number of non-overlapping partitions of $[n]$; this number is
called the $n$th \emph{Bessel number} \cite[p. 423]{FlSc90}.  The
first few Bessel numbers are
$$
1,1,2,5,14,43,143,509,1922,7651,31965,139685,636712. $$
We denote
by $S^*(n,k)$ the number of non-overlapping partitions of $[n]$ into
$k$ blocks.

\subsection*{Involutions}

An \emph{involution} is a permutation which is its own inverse. We
denote by $I_n$ the number of involutions in $\sym_n$. The sequence
$\{I_n\}_0^{\infty}$ starts with
$$
1,1,2,4,10,26,76,232,764,2620,9496,35696,140152. $$

\subsection*{Dyck paths}

A \emph{Dyck path} of length $2n$ is a lattice path from $(0,0)$ to
$(2n,0)$ with steps $(1,1)$ and $(1,-1)$ that never goes below the
$x$-axis. Letting $u$ and $d$ represent the steps $(1,1)$ and $(1,-1)$
respectively, we code such a path with a word over $\{u,d\}$. For
example, the path
\begin{center}
  \includegraphics[width=25ex]{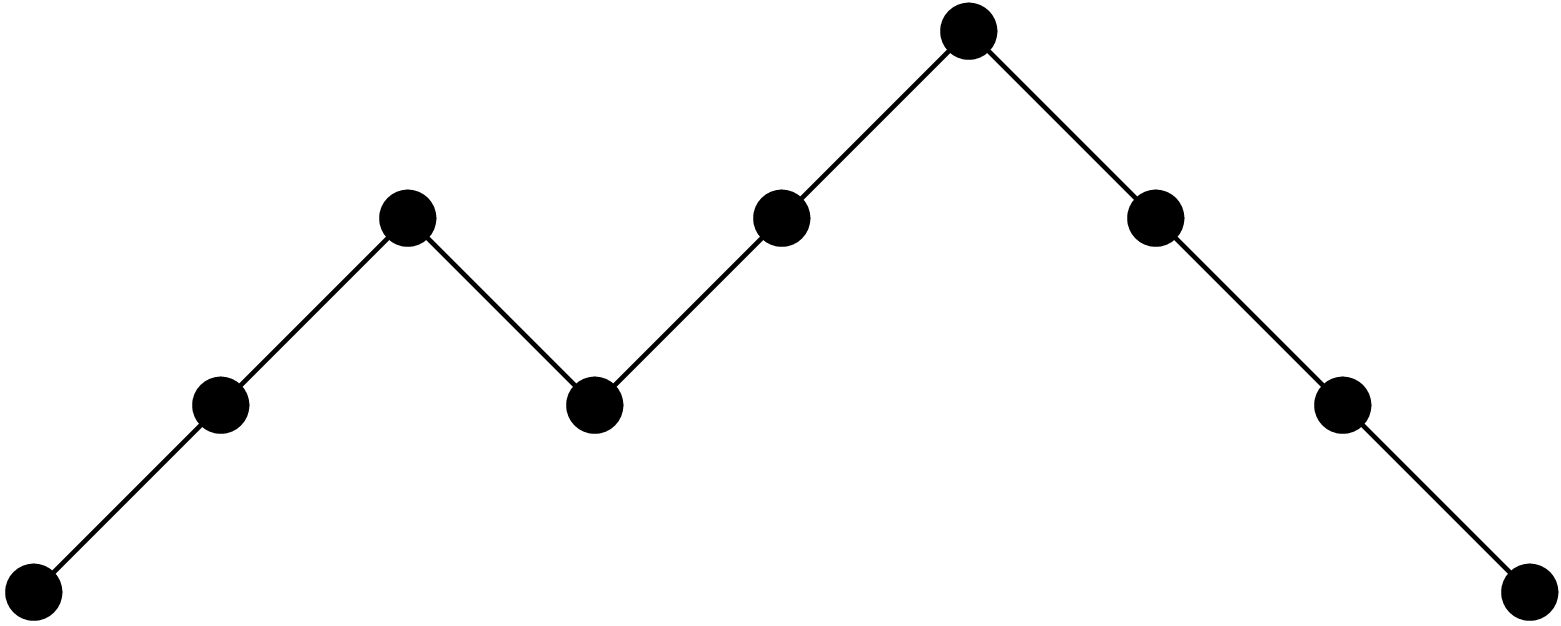}
\end{center}
is coded by $uuduuddd$. The $n$th \emph{Catalan number} $C_n=\frac 1
{n+1} \binom {2n} n$ counts the number of Dyck paths of length $2n$.
The sequence of Catalan numbers starts with
$$1,1,2,5,14,42,132,429,1430,4862,16796,58786,208012.$$

\subsection*{Motzkin paths}

A \emph{Motzkin path} of length $n$ is a lattice path from $(0,0)$ to
$(n,0)$ with steps $(1,0)$, $(1,1)$, and $(1,-1)$ that never goes
below the $x$-axis. Letting $\ell$, $u$, and $d$ represent the steps
$(1,0)$, $(1,1)$, and $(1,-1)$ respectively, we code such a path with
a word over $\{\ell, u,d\}$. For example, the path
\begin{center}
  \includegraphics[width=25ex]{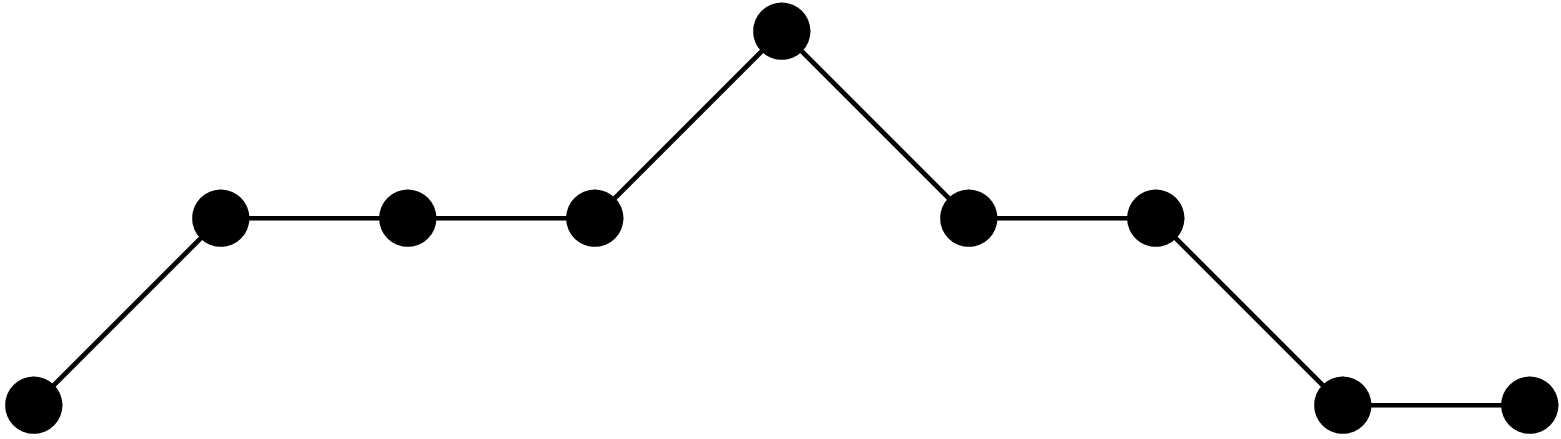}
\end{center}
is coded by $u\ell \ell ud\ell d\ell $. The $n$th \emph{Motzkin
  number} $M_n$ is the number of Motzkin paths of length $n$. The
first few of the Motzkin numbers are
$$
1,1,2,4,9,21,51,127,323,835,2188,5798,15511. $$

\section{Three classes of patterns}

Let $\pi=a_1 a_2 \cdots a_n \in \sym_n$. Define the \emph{reverse} of
$\pi$ as $\pi^r:=a_n \cdots a_2 a_1$, and define the \emph{complement}
of $\pi$ by $\pi^c(i) = n+1-\pi(i)$, where $i\in[n]$.

\begin{proposition} 
  With respect to being equidistributed, the twelve pattern statistics
  of length three with one dash fall into the following three classes.
  \begin{itemize}
  \item[(i)] \hspace{.8ex}\axbc,\; \cxba,\; \abxc,\; \cbxa.
  \item[(ii)] \hspace{.8ex}\axcb,\; \cxab,\; \baxc,\; \bcxa.
  \item[(iii)] \hspace{.8ex}\bxac,\; \bxca,\; \acxb,\; \caxb.
  \end{itemize}
\end{proposition}

\begin{proof}
  The bijections $\pi\mapsto\pi^r$, $\pi\mapsto\pi^c$, and
  $\pi\mapsto(\pi^r)^c$ give the equidistribution part of the result.
  Calculations show that these three distributions differ pairwise on
  $\sym_4$.
\end{proof}

% ------------------------------------------------------------------------
\section{Permutations avoiding a pattern of class one or two}

\begin{proposition}\label{bell1}
  Partitions of $[n]$ are in one-to-one correspondence with
  $(\axbc)$-avoiding permutations in $\sym_n$. Hence
  $|\sym_n(\axbc)|=B_n$.
\end{proposition}

\begin{proof}[First proof]
  Recall that the Bell numbers satisfy $B_0=1$, and
  $$B_{n+1}= \sum_{k=0}^n \binom n k B_k.$$
  We show that
  $|\sym_n(\axbc)|$ satisfy the same recursion. Clearly,
  $\sym_0(\axbc) = \{\emptyword\}$.  For $n > 0$, let
  $M=\{2,3,\ldots,n+1\}$, and let $S$ be a $k$ element subset of $M$.
  For each $(\axbc)$-avoiding permutation $\sigma$ of $S$ we construct
  a unique $(\axbc)$-avoiding permutation $\pi$ of $[n+1]$. Let $\tau$
  be the word obtained by writing the elements of $M\setminus S$ in
  decreasing order. Define $\pi := \sigma 1 \tau$.
  
  Conversely, if $\pi = \sigma 1 \tau$ is a given $(\axbc)$-avoiding
  permutation of $[n+1]$, where $|\sigma|=k$, then the letters of
  $\tau$ are in decreasing order, and $\sigma$ is an $(\axbc)$-avoiding
  permutation of the $k$ element set $\{2,3,\ldots,n+1\}\setminus \{i
  : i \text{ is a letter in }\tau\}$.
\end{proof}

\begin{proof}[Second proof]
  Given a partition $\pi$ of $[n]$, we introduce a standard
  representation of $\pi$ by requiring that:
  \begin{itemize}
  \item[(a)] Each block is written with its least element first, and
    the rest of the elements of that block are written in decreasing
    order.
  \item[(b)] The blocks are written in decreasing order of their least
    element, and with dashes separating the blocks.
  \end{itemize}
  
  Define $\wpi$ to be the permutation we obtain from $\pi$ by writing
  it in standard form and erasing the dashes. We now argue that
  $\wpi:=a_1a_2\cdots a_n$ avoids $(\axbc)$.  If $a_i<a_{i+1}$, then
  $a_i$ and $a_{i+1}$ are the first and the second element of some
  block. By the construction of $\wpi$, $a_i$ is a left-to-right
  minimum, hence there is no $j\in[i-1]$ such that $a_j<a_i$.
  
  Conversely, $\pi$ can be recovered uniquely from $\wpi$ by inserting
  a dash in $\wpi$ preceding each left-to-right minimum, apart from
  the first letter in $\wpi$.  Thus $\pi\mapsto\wpi$ gives the desired
  bijection.
\end{proof}

\begin{example}
  As an illustration of the map defined in the above proof, let
  $$\pi=\{\{1,3,5\},\{2,6,9\},\{4,7\},\{8\}\}.$$
  Its standard form is
  $8 \dd 47 \dd 296 \dd 153$. Thus $\wpi= 847296153$.
\end{example}

\begin{porism}
  Let $L(\pi)$ be the number of left-to-right minima of $\pi$. Then
  $$
  \sum_{\pi\in\sym_n(\axbc)} \!\!\!\!\!\!\!\! x^{L(\pi)} =
  \sum_{k\geq 0} S(n,k)x^k.
  $$
\end{porism}

\begin{proof}
  This result follows readily from the second proof of
  Proposition~\ref{bell1}. We here give a different proof, which is
  based on the fact that the Stirling numbers of the second kind
  satisfy
  \begin{equation*}
    S(n,k) = S(n-1,k-1) + kS(n-1,k).
  \end{equation*}
  
  Let $T(n,k)$ be the number of permutations in $\sym_n(\axbc)$ with
  $k$ left-to-right minima. We show that the $T(n,k)$ satisfy the same
  recursion as the $S(n,k)$.
  
  Let $\pi$ be an $(\axbc)$-avoiding permutation of $[n-1]$. To insert
  $n$ in $\pi$, preserving $(\axbc)$-avoidance, we can put $n$ in
  front of $\pi$ or we can insert $n$ immediately after each
  left-to-right minimum.  Putting $n$ in front of $\pi$ creates a new
  left-to-right minimum, while inserting $n$ immediately after a
  left-to-right minimum does not.
\end{proof}

% ------------------------------------------------------------------------

\begin{proposition}\label{bell2}
  Partitions of $[n]$ are in one-to-one correspondence with
  $(\axcb)$-avoiding permutations in $\sym_n$. Hence
  $|\sym_n(\axcb)|=B_n$.
\end{proposition}

\begin{proof}
  Let $\pi$ be a partition of $[n]$. We introduce a standard
  representation of $\pi$ by requiring that:
  \begin{itemize}
  \item[(a)] The elements of a block are written in increasing order.
  \item[(b)] The blocks are written in decreasing order of their least
    element, and with dashes separating the blocks.
  \end{itemize}
  Notice that this standard representation is different from the one
  given in the second proof of Proposition~\ref{bell1}. Define $\wpi$
  to be the permutation we obtain from $\pi$ by writing it in standard
  form and erasing the dashes. It easy to see that $\wpi$ avoids
  $(\axcb)$.  Conversely, $\pi$ can be recovered uniquely from $\wpi$
  by inserting a dash in between each descent in $\wpi$.
\end{proof}

\begin{example}
  As an illustration of the map defined in the above proof, let
  $$\pi=\{\{1,3,5\},\{2,6,9\},\{4,7\},\{8\}\}.$$
  Its standard form is
  $8 \dd 47 \dd 269 \dd 135$. Thus $\wpi= 847269135$.
\end{example}

\begin{porism}
  $$
  \sum_{\pi\in\sym_n(\axcb)} \!\!\!\!\!\!\!\! x^{1+\des \pi } =
  \sum_{k\geq 0} S(n,k)x^k.
  $$
\end{porism}

\begin{proof}
  From the proof of Proposition~\ref{bell2} we see that $\pi$ has
  $k+1$ blocks precisely when $\wpi$ has $k$ descents.
\end{proof}

% ------------------------------------------------------------------------

\begin{proposition}\label{involutions}
  Involutions in $\sym_n$ are in one-to-one correspondence with
  permutations in $\sym_n$ that avoid $(\axbc)$ and $(\axcb)$.  Hence
  $$ |\sym_n(\axbc, \axcb)|=I_n.$$
\end{proposition}

\begin{proof}
  We give a combinatorial proof using a bijection that is essentially
  identical to the one given in the second proof of
  Proposition~\ref{bell1}.
  
  Let $\pi\in\sym_n$ be an involution. Recall that $\pi$ is an
  involution if and only if each cycle of $\pi$ is of length one or
  two. We now introduce a standard form for writing $\pi$ in cycle
  notation by requiring that:
  \begin{itemize}
  \item[(a)] Each cycle is written with its least element first.
  \item[(b)] The cycles are written in decreasing order of their least
    element.
  \end{itemize}
  Define $\wpi$ to be the permutation obtained from $\pi$ by writing
  it in standard form and erasing the parentheses separating the
  cycles.

  Observe that $\wpi$ avoids $(\axbc)$: Assume that $a_i<a_{i+1}$,
  that is $(a_i\,a_{i+1})$ is a cycle in $\pi$, then $a_i$ is a
  left-to-right minimum in $\pi$. This is guaranteed by the
  construction of $\wpi$. Thus there is no $j<i$ such that $a_j<a_i$.
  
  The permutation $\wpi$ also avoids $(\axcb)$: Assume that
  $a_i>a_{i+1}$, then $a_{i+1}$ must be the smallest element of some
  cycle. Then $a_i$ is a left-to-right minimum in $\pi$.

  Conversely, if $\wpi:=a_1\ldots a_n$ is an $\{\axbc, \axcb\}$-avoiding
  permutation then the involution $\pi$ is given by: $(a_i\,a_{i+1})$ is a
  cycle in $\pi$ if and only if $a_i<a_{i+1}$.
\end{proof}

\begin{example}
  The involution $\pi=826543719$ written in standard form is
  $$(9)(7)(4\,5)(3\,6)(2)(1\,8),$$
  and hence $\wpi = 974536218$.
\end{example}

\begin{porism}
  $$\sym_n(\axbc, \axcb) = \sym_n(\axbc, acb)
                         = \sym_n(abc, \axcb)
                         = \sym_n(abc, acb).$$
\end{porism}

\begin{proof}
  Kitaev \cite{Ki00} observed that the dashes in the patterns $(\axbc)$
  and $(\axcb)$ are immaterial for the proof of
  Proposition~\ref{involutions}. The result may, however, also be proved
  directly. For an example of such a proof see the proof of
  Lemma~\ref{lemma2}.
\end{proof}

\begin{porism}\label{des_fix}
  The number of permutations in $\sym_{n+k}(\axbc,\axcb)$ with $n-1$
  descents equals the number of involutions in $\sym_{n+k}$ with $n-k$
  fixed points.
\end{porism}

\begin{proof}
  Under the bijection $\pi\mapsto\wpi$ in the proof of
  Proposition~\ref{involutions}, a cycle of length two in $\pi$
  corresponds to an occurrence of $(\ab)$ in $\wpi$. Hence, if $\pi$
  has $n-2k$ fixed points, then $\wpi$ has $n-k-1$ descents.
  Substituting $n+k$ for $n$ we get the desired result.
\end{proof}

To take the analysis of descents in $\{\axbc,\axcb\}$-avoiding
permutations further, we introduce the polynomial
$$
A_n(x) = \!\!\! \!\!\!  \sum_{\pi\in\sym_n(\axbc,\axcb)}
\!\!\!\!\!\! \!\!\!\!\!\! \!  x^{1+\des\,\pi},
$$
and call it the $n$th \emph{Eulerian polynomial for
  $\{\axbc,\axcb\}$-avoiding permutations}. Direct enumeration shows
that the sequence $\{A_n(x)\}$ starts with
$$
\begin{array}{rclclclclclclclcl}
  A_0(x) &=& 1 \\
  A_1(x) &=& &\!\!+\!\!& x\\
  A_2(x) &=& & & x &\!\!+\!\!& x^2\\
  A_3(x) &=& & & & & 3x^2 &\!\!+\!\!& x^3\\
  A_4(x) &=& & & & & 3x^2 &\!\!+\!\!& 6x^3  &\!\!+\!\!& x^4\\
  A_5(x) &=& & & & & & & 15x^3 &\!\!+\!\!& 10x^4 &\!\!+\!\!& x^5\\
  A_6(x) &=& & & & & & & 15x^3 &\!\!+\!\!& 45x^4 &\!\!+\!\!& 15x^5 &\!\!+\!\!& x^6\\
  A_7(x) &=& & & & & & & & & 105x^4 &\!\!+\!\!& 105x^5 &\!\!+\!\!& 21x^6 &\!\!+\!\!& x^7.\\
\end{array}
$$
We will relate these polynomials to the so called Bessel
polynomials.  The $n$th \emph{Bessel polynomial} $y_n(x)$ is defined
by
\begin{equation}\label{bessel_expl}
  y_n(x) = \sum_{k=0}^n \binom {n+k} k \binom n k \frac {k!} {2^k}
  x^k.
\end{equation}
The first six of the Bessel polynomials are
$$
\begin{array}{rclclclclclclcl}
  y_0(x) &=& 1 \\
  y_1(x) &=& 1 &\!\!+\!\!& x\\
  y_2(x) &=& 1 &\!\!+\!\!& 3x  &\!\!+\!\!& 3x^2\\
  y_3(x) &=& 1 &\!\!+\!\!& 6x  &\!\!+\!\!& 15x^2  &\!\!+\!\!& 15x^3\\
  y_4(x) &=& 1 &\!\!+\!\!& 10x &\!\!+\!\!& 45x^2  &\!\!+\!\!& 105x^3 &\!\!+\!\!& 105x^4\\
  y_5(x) &=& 1 &\!\!+\!\!& 15x &\!\!+\!\!& 105x^2 &\!\!+\!\!& 420x^3 &\!\!+\!\!& 945x^4 &\!\!+\!\!& 945x^5.\\
\end{array}
$$
These polynomials satisfy the second order differential
equation
$$x^2\frac {d^2y} {dx^2} + 2(x+1) \frac {dy} {dx} = n(n+1)y.$$
Moreover, the Bessel polynomials satisfy the recurrence relation
\begin{equation}\label{bessel_rec}
  y_{n+1}(x) = (2n+1)xy_n(x) + y_{n-1}(x).
\end{equation}

\begin{proposition}\label{euler}
  Let $y_n(x)$ be the $n$th Bessel polynomial, and let $A_n(x)$ be the
  $n$th Eulerian polynomial for $\{\axbc,\axcb\}$-avoiding
  permutations. Then
  \begin{itemize}
  \item[(i)] $\sum_n y_n(x) (xt)^n$ generates $\{A_n(t)\}$, that is
    $$\sum_{n\geq 0} A_n(t) x^n = \sum_{n\geq 0} y_n(x) (xt)^n.$$
  \item[(ii)] $A_0(x) = 1$, $A_1(x) = x$, and for $n\geq2$, we have
    $$A_{n+2}(x) = x(1+x+2x\frac d {dx}) A_n(x).$$
  \item[(iii)] $A_n(x)$ is explicitly given by
    $$A_n(x) = \sum_{k=0}^n \binom n k \binom {n-k} k \frac {k!} {2^k}
    x^{n-k}.$$
  \end{itemize}
\end{proposition}

\begin{proof}
  Let $I_n^k$ denote the number of involutions in $\sym_n$ with $k$
  fixed points. Then Porism~\ref{des_fix} is equivalently stated as
  \begin{equation}\label{eulerian_poly}
    A_n(x) = \sum_{k\geq 0} I_n^{2k-n} x^k.
  \end{equation}
  In \cite{DuFa91} Dulucq and Favreau showed that the Bessel
  polynomials are given by
  \begin{equation}\label{bessel_poly}
    y_n(x) = \sum_{k\geq 0} I_{n+k}^{n-k} x^k.
  \end{equation}
  To prove (i), multiply Equation~(\ref{bessel_poly}) by $(xt)^n$ and
  sum over $n$.
  \begin{align*}
    \sum_{n\geq 0} y_n(x)(xt)^n
    &\,=\, \sum_{n\geq 0} \sum_{k\geq 0} I_{n+k}^{n-k} t^n x^{n+k}  \\
    &\,=\, \sum_{k\geq 0} \sum_{n\geq 0} I_{k}^{2n-k} t^n x^k
    && \text{By substituting $n-k$ for $k$.}\\
    &\,=\, \sum_{k\geq 0} A_k(t) x^k
    && \text{By Equation~(\ref{eulerian_poly}).}
  \end{align*}
  We now multiply Equation~(\ref{bessel_rec}) by $(xt)^n$ and sum over
  $n$. Tedious but straightforward calculations then yield (ii) from
  (i). Finally, we obtain (iii) from Equation~(\ref{bessel_expl}) by
  identifying coefficients in (i).
\end{proof}

% ------------------------------------------------------------------------

\begin{definition}
  Let $\pi$ be an arbitrary partition whose non-singleton blocks
  $\{A_1,\ldots,A_k\}$ are ordered so that for all $i\in[k-1]$, $\min
  A_i>\min A_{i+1}$. If $\max A_i>\max A_{i+1}$ for all $i\in [k-1]$,
  then we call $\pi$ a \emph{monotone partition}. The set of monotone
  partitions of $[n]$ is denoted by $\mono_n$.
\end{definition}

\begin{example}
  The partition
  $$\pi \,=\, \FIGmono $$
  is monotone.
\end{example}

\begin{proposition}\label{mono}
  Monotone partitions of $[n]$ are in one-to-one correspondence with
  permutations in $\sym_n$ that avoid $(\axbc)$ and $(\abxc)$. Hence
  $$|\sym_n(\axbc,\abxc)| = |\mono_n|.$$
\end{proposition}

\begin{proof} 
  Given $\pi$ in $\mono_n$, let $A_1\dd A_2\dd\cdots\dd A_k$ be the
  result of writing $\pi$ in the standard form given in the second
  proof of Proposition~\ref{bell1}, and let $\wpi=A_1 A_2\cdots A_k$.
  By the construction of $\wpi$ the fist letter in each $A_i$ is a
  left-to-right minimum. Furthermore, since $\pi$ is monotone the
  second letter in each non-singleton $A_i$ is a right-to-left
  maximum.  Therefore, if $xy$ is an $(ab)$-subword of $\wpi$, then $x$
  is left-to-right minimum and $y$ is a right-to-left maximum. Thus
  $\wpi$ avoids both $(\axbc)$ and $(\abxc)$.
  
  Conversely, given $\wpi$ in $\sym_n(\axbc,\abxc)$, let $A_1\dd
  A_2\dd\cdots\dd A_k$ be the result of inserting a dash in $\wpi$
  preceding each left-to-right minimum, apart from the first letter in
  $\wpi$. Since $\wpi$ is $(\abxc)$-avoiding, the second letter in
  each non-singleton $A_i$ is a right-to-left maximum. The second
  letter in $A_i$ is the maximal element of $A_i$ when $A_i$ is viewed
  as a set. Thus $\pi = \{A_1,A_2,\ldots,A_k \}$ is monotone.
\end{proof}

What is left for us to show is that there is a one-to-one
correspondence between monotone partitions and non-overlapping
partitions. The proof we give is strongly influenced by the work of
Flajolet \cite{FlSc90}.

\begin{proposition}\label{nop-mono}
  Monotone partitions of $[n]$ are in one-to-one correspondence with
  non-overlapping partitions of $[n]$. Hence $|\mono_n| = B_n^*$.
\end{proposition}

\begin{proof}
  If $k$ is the minimal element of a non-singleton block, then call
  $k$ the \emph{first} element of that block. Similarly, If $k$ is the
  maximal element of a non-singleton block, then call $k$ the
  \emph{last} element of that block. An element of a non-singleton
  that is not a first or last element is called an \emph{intermediate}
  element. Let us introduce an ordering of the blocks of a partition. A
  block $A$ is smaller than a block $B$ if $\min A < \min B$.
  
  We define a map $\Phi$ that to each non-overlapping partition $\pi$
  of $[n]$ gives a unique monotone partition $\Phi(\pi)$ of
  $[n]$.  Let the integer $k$ range from $1$ to $n$.
  \begin{itemize}
  \item[(a)] If $k$ is the first element of a block of $\pi$, then
    open a new block in $\Phi(\pi)$ by letting $k$ be its first
    element. (A block $B$ is open if $\max B<k$.)
  \item[(b)] If $k$ is the last element of a block of $\pi$, then
    close the smallest open block of $\Phi(\pi)$ by letting $k$ be its
    last element.
  \item[(c)] If $k$ is an intermediate element of some block $B$ of
    $\pi$, and $B$ is the $i$:th largest open block of $\pi$, then let
    $k$ belong to the $i$th largest open block of $\Phi(\pi)$.
  \item[(d)] If $\{k\}$ is a singleton block of $\pi$, then let
    $\{k\}$ be a singleton block of $\Phi(\pi)$.
  \end{itemize}

  Observe that $\Phi(\pi)$ is monotone. Indeed, it is only in (b) that
  we close a block of $\Phi(\pi)$, and we always close the smallest
  open block of $\Phi(\pi)$.

  Conversely, we give a map $\Psi$ that to each monotone partition
  $\pi$ be a of $[n]$ gives a unique non-overlapping partition
  $\Psi(\pi)$ of $[n]$. Define $\Psi$ the same way as $\Phi$ is
  defined, except for case (c), where we instead of closing the
  smallest open block close the largest open block.

  It is easy to see that $\Phi$ and $\Psi$ are each others inverses
  and hence they are bijections.
\end{proof}

\begin{corollary}\label{nop}
  The NOPs (non-overlapping partitions) of $[n]$ are in one-to-one
  correspondence with permutations in $\sym_n$ that avoid $(\axbc)$
  and $(\abxc)$. Hence
  $$|\sym_n(\axbc, \abxc)|=B^*_n.$$
\end{corollary}

\begin{proof}
  Follows immediately from Proposition~\ref{mono} together with
  Proposition~\ref{nop-mono}.
\end{proof}

\begin{example}
  By the proof of Proposition~\ref{nop-mono}, the non-overlapping partition
  \begin{align*}
    \pi       &\,=\, \FIGnop \\
    \intertext{corresponds to the monotone partition}
    \Phi(\pi) &\,=\, \FIGmono \\
    \intertext{that according to the proof of Proposition~\ref{mono} 
      corresponds to the $\{\axbc,\abxc\}$-avoiding permutation}
    \widehat{\Phi(\pi)} &\,=\, 1\ns 0\;1\ns 3\;1\ns 1\;9\;4\;1\ns 2\;6\;3\;8\;1\;7\;5\;2.
  \end{align*}
\end{example}

\begin{porism}
  Let $L(\pi)$ be the number of left-to-right minima of $\pi$. Then
  $$
  \sum_{\pi\in\sym_n(\axbc,\abxc)} \!\!\!\!\!\!\!\!\!\!\!\!\!\!
  x^{L(\pi)} = \sum_{k\geq 0} S^*(n,k)x^k.
  $$
\end{porism}

\begin{proof}
  Under the bijection $\pi\mapsto\wpi$ in the proof of
  Proposition~\ref{mono}, the number of blocks in $\pi$ determines the
  number of left-to-right minima of $\wpi$, and vice versa.  The
  number of blocks is not changed by the bijection $\Phi_1\circ\Psi_2$
  in the proof of Proposition~\ref{nop-mono}.
\end{proof}

% ------------------------------------------------------------------------
\section{Permutations avoiding a pattern of class three}

In \cite{Kn73v1} Knuth observed that there is a
one-to-one correspondence between $(\bxaxc)$-avoiding permutations and
Dyck paths. For completeness and future reference we give this result
as a lemma, and prove it using one of the least known bijections.
First we need a definition. For each word $x=x_1x_2\cdots x_n$
without repeated letters, we define the \emph{projection} of $x$ onto
$S_n$, which we denote $\proj(x)$, by
$$
\proj(x) = a_1 a_2 \cdots a_n \,,\;\text{ where }\; a_i = |\{j\in
[n] : x_i\geq x_j \}|.
$$
Equivalently, $\proj(x)$ is the permutation in $\sym_n$ which is
order equivalent to $x$. For example, $\proj(265) = 132$.

\begin{lemma}\label{lemma1}
  $|\sym_n(\bxaxc)|=C_n.$
\end{lemma}

\begin{proof}
  Let $\pi=a_1 a_2 \cdots a_n$ be a permutation of $[n]$ such that
  $a_k=1$. Then $\pi$ is $(\bxaxc)$-avoiding if and only if
  $\pi=\sigma 1 \tau$, where $\sigma:=a_1\cdots a_{k-1}$ is a
  $(\bxaxc)$-avoiding permutation of $\{n,n-1,\ldots,n-k+1\}$, and
  $\tau:=a_{k+1}\cdots a_n$ is a $(\bxaxc)$-avoiding permutation of
  $\{2,3,\ldots,k\}$.
  
  We define recursively a mapping $\Phi$ from $\sym_n(\bxaxc)$ onto
  the set of Dyck paths of length $2n$. If $\pi$ is the empty word,
  then so is the Dyck path determined by $\pi$, that is,
  $\Phi(\emptyword) = \emptyword$. If $\pi \neq \emptyword$, then we
  can use the factorisation $\pi=\sigma 1 \tau$ from above, and define
  $\Phi(\pi) = u\,(\Phi\circ\proj)(\sigma)\,d\,(\Phi\circ\proj)(\tau)$. It
  is easy to see that $\Phi$ may be inverted, and hence is a
  bijection.
\end{proof}

\begin{lemma}\label{lemma2}
  A permutation avoids $(\bxac)$ if and only if it avoids $(\bxaxc)$.
\end{lemma}

\begin{proof} 
  The sufficiency part of the proposition is trivial. The necessity
  part is not difficult either. Assume that $\pi$ contains a
  $(\bxaxc)$-subword. Then there exist
  $$
  A,B,C, n_1,n_2,\ldots,n_r \in [n], \text{ where } A < B < C,
  $$
  such that $B\hspace{-1pt}AC$ is a subword of $\pi$, and $A
  n_1\cdots n_r C$ is a segment of $\pi$. If $n_1>B$, then
  $B\hspace{-1pt}A n_1$ form a $(\bxac)$-subword in $\pi$. Assume that
  $n_1<B$.  Indeed, to avoid forming a $(\bxac)$-subword we will have
  to assume that $n_i<B$ for all $i\in[r]$, but then $B n_r C$ is a
  $(\bxac)$-subword. Accordingly we conclude that there exists at
  least one $(\bxac)$-subword in $\pi$.
\end{proof}

\begin{proposition}\label{catalan}
  Dyck paths of length $2n$ are in one-to-one correspondence with
  $(\bxac)$-avoiding permutations in $\sym_n$. Hence
  $$|\sym_n(\bxac)|=\frac 1 {n+1} \binom {2n} n.$$
\end{proposition}

\begin{proof} 
  Follows immediately from Lemma~\ref{lemma1} and Lemma~\ref{lemma2}.
\end{proof}

\begin{proposition}
  Let $L(\pi)$ be the number of left-to-right minima of $\pi$. Then
  $$
  \sum_{\pi\in\sym_n(\bxac)} \!\!\!\!\!\!\!\! x^{L(\pi)} =
  \sum_{k\geq 0} \frac k {2n-k} \binom{2n-k} n x^k.
  $$
\end{proposition}

\begin{proof}
  A \emph{return step} in a Dyck path $\delta$ is a $d$ such that
  $\delta = \alpha u \beta d \gamma$, for some Dyck paths $\alpha$,
  $\beta$, and $\gamma$. A useful observation is that every non-empty
  Dyck path $\delta$ can be uniquely decomposed as $\delta = u \alpha
  d \beta$, where $\alpha$ and $\beta$ are Dyck paths.  This is the
  so-called \emph{first return decomposition} of $\delta$. Let
  $R(\delta)$ denote the number of return steps in $\delta$. 
  
  In \cite{De99} Deutsch showed that the distribution of $R$ over all
  Dyck paths of length $2n$ is the distribution we claim that $L$ has
  over $\sym_n(\bxac)$.
  
  Let $\gamma$ be a Dyck path of length $2n$, and let $\gamma = u
  \alpha d \beta$ be its first return decomposition. Then $R(\gamma) =
  1 + R(\beta)$.  Let $\pi\in\sym_n(\bxac)$, and let $\pi=\sigma 1
  \tau$ be the decomposition given in the proof of Lemma~\ref{lemma1}.
  Then $L(\pi) = 1 + L(\sigma)$.  The result now follows by induction.
\end{proof}

In addition, it is easy to deduce that left-to-right minima,
left-to-right maxima, right-to-left minima, and right-to-left maxima
all share the same distribution over $\sym_n(\bxac)$.

% ------------------------------------------------------------------------

\begin{proposition}
  Motzkin paths of length $n$ are in one-to-one correspondence with
  permutations in $\sym_n$ that avoid $(\axbc)$ and $(\acxb)$. Hence
  $$|\sym_n(\axbc, \acxb)|=M_n.$$
\end{proposition}

\begin{proof}
  We mimic the proof of Lemma~\ref{lemma1}. Let $\pi\in\sym_n(\axbc,
  \acxb)$.  Since $\pi$ avoids $(\acxb)$ it also avoids $(\axcxb)$ by
  Lemma~\ref{lemma2} via $\pi\mapsto (\pi^c)^r$. Thus we may write
  $\pi = \sigma n \tau$, where $\pi(k)=n$, $\tau$ is an $\{\axbc,
  \acxb\}$-avoiding permutation of $\{n-1,n-2,\ldots,n-k+1\}$, and
  $\tau$ is an $\{\axbc, \acxb\}$-avoiding permutation of $[n-k]$.  If
  $\sigma\neq\emptyword$ then $\sigma=\sigma'r$ where $r=n-k+1$, or
  else an $(\axbc)$-subword would be formed with $n$ as the '$c$' in
  $(\axbc)$. Define a map $\Phi$ from $\sym_n(\axbc, \acxb)$ to the
  set of Motzkin paths by $\Phi(\emptyword) = \emptyword$ and
  $$\Phi(\pi) =
  \begin{cases}
    \ell\,(\Phi\circ\proj)(\sigma) & \text{ if } \pi = n\sigma, \\
    u\,(\Phi\circ\proj)(\sigma)\,d\,\Phi(\tau) & \text{ if }
                         \pi = \sigma r n \tau \text{ and } r=n-k+1.\\
  \end{cases}
  $$
  Its routine to find the inverse of $\Phi$.
\end{proof}

\begin{example}
  Let us find the Motzkin path associated to the 
  $\{\axbc,\acxb\}$-avoiding permutation $76453281$.
  \begin{eqnarray*}
    \Phi(76453{\bf 28}1) &=& u \Phi({\bf 5}4231) d \Phi({\bf 1}) \\
                         &=& u\ell \Phi({\bf 4}231) d\ell \\
                         &=& u\ell\ell \Phi({\bf 23}1) d\ell \\
                         &=& u\ell\ell u d \Phi({\bf 1}) d\ell \\
                         &=& u\ell\ell u d\ell d\ell \\
  \end{eqnarray*}
\end{example}

\section*{Acknowledgement}

I am greatly indebted to my advisor Einar Steingr\'{\i}msson, who put
his trust in me and gave me the opportunity to study mathematics on a
postgraduate level. This work has benefited from his knowledge,
enthusiasm and generosity.

\bibliographystyle{plain}
\nocite{SlPl95}
\nocite{St97}
\bibliography{gpa}

\end{document}